\documentclass[11pt]{amsart}
\usepackage{amssymb, amsthm, amsmath, gensymb}
\usepackage{graphicx, comment}
\usepackage[small]{caption}
\usepackage{subcaption}
\usepackage{epsfig}
\usepackage{tikz, pgfplots, float}

\newtheorem{thm}{Theorem}[section]

\theoremstyle{definition}

\newcommand{\F}{\mathcal F}
\newcommand{\HH}{\mathcal H}
\newcommand{\G}{\mathcal G}
\newcommand{\bC}{\mathbf C}
\newcommand{\C}{\mathcal C}
\newcommand{\PP}{\mathcal P}

\newcommand{\LAw}{\mathop{}\!\mathrm{La}}
\newcommand{\La}{\mathop{}\!\mathrm{La}^*}

\title{Non-principality of forbidden subposet problems}
\author{Bal\'azs Patk\'os}
\date{}

\begin{document}

\begin{abstract}
 We show pairs of posets such that the largest families not containing any of them is substantially smaller than the largest families avoiding only one of them.  
\end{abstract}

\maketitle

\section{Introduction}

We use standard notation: $[n]$ denotes the set of the first $n$ positive integers, and for a set $X$ we denote by $2^X$ and $\binom{X}{k}$ the power set of $X$ and the family of all subsets of $X$ of size $k$, respectively.

This note is about the area of forbidden subposet problems. For a more detailed introduction than what we have below, we refer the reader to the surveys \cite{AMP,GLsurv} and chapter 7 of \cite{GP}. A family $\G$ of sets is a \textit{weak copy} of a poset $P=(P,\leqslant)$ if there exists a bijection $\iota: P\rightarrow \G$ such that $p\leqslant q$ implies $\iota(p)\subseteq \iota(q)$. If the bijection satisfies the property $p\leqslant q$ if and only if $\iota(p)\subseteq \iota(q)$, then $\G$ is a \textit{strong copy} of $P$. A family $\F$ is \textit{$P$-free} if it contains no weak copy of $P$, and $\F$ is \textit{strong $P$-free} if it contains no strong copy of $P$. The maximum size of a $P$-free / strong $P$-free family $\F\subseteq 2^{[n]}$ is denoted by $\LAw(n,P) / \La(n,P)$. The problem  of determining these quantities was introduced by Katona and Tarj\'an \cite{KatonaTarjan}. These notions and parameters can be introduced for families $\PP$ of posets: $\F$ is (strong) $\PP$-free if it is (strong) $P$-free for all $P\in \PP$. The order of magnitude of $\LAw(n,P)$ and $\La(n,P)$ is $\Theta(\binom{n}{\lfloor \frac{n}{2}\rfloor})$ \cite{MP} provided $P$ contains at least one pair of elements in relation, but it is not known if $\pi(P)=\lim_{n\rightarrow \infty}\frac{\LAw(n,P)}{\binom{n}{\lfloor \frac{n}{2}\rfloor}}$ and $\pi^*(P)=\lim_{n\rightarrow \infty}\frac{\La(n,P)}{\binom{n}{\lfloor \frac{n}{2}\rfloor}}$ exist. We will use the notation $\pi^-(P)=\liminf_{n\rightarrow \infty}\frac{\LAw(n,P)}{\binom{n}{\lfloor \frac{n}{2}\rfloor}} $ and $\pi^{*,-}(P)=\liminf_{n\rightarrow \infty}\frac{\La(n,P)}{\binom{n}{\lfloor \frac{n}{2}\rfloor}} $.

An easy way to construct $P$-free families is to consider as many consecutive middle levels of $2^{[n]}$ as possible without creating a copy of $P$. Formally, we define $e(P)$ and $e^*(P)$ to be the largest integers $\ell$ such that $\cup_{i=j}^{j+\ell-1}\binom{[n]}{i}$ is $P$-free (strong $P$-free) for all $j$ and $n$. By definition and by the fact that $\binom{n}{\lfloor\frac{n}{2}\rfloor+j}=(1+o(1))\binom{n}{\lfloor\frac{n}{2}\rfloor}$ for any fixed $j$, we have $\pi^-(P)\ge e(P)$ and $\pi^{*,-}(P)\ge e^*(P)$. These definitions and inequalities go through for families $\PP$ of posets, too. Ellis, Ivan, and Leader \cite{EIL} were the first to prove that there exist posets for which the inequalities are strict. Recently, there has been substantial progress \cite{P,PT1,PT2,T} in the area. 

Principality and non-principality are phenomena observed in other areas of extremal combinatorics. A consequence of the Erd\H os-Stone-Simonovits theorem and the Tur\'an graph construction is that the Tur\'an density of a family $\HH$ of forbidden subgraphs is the minimum of the Tur\'an densities of the graphs $H\in \HH$. This does not remain valid if one considers $r$-graphs with $r\ge 3$ as shown by Balogh \cite{Balogh2002} for the case $r=3$ and by Mubayi and Pikhurko \cite{MubayiPikhurko2008} for general $r$. Forbidden subposet problems more closely resemble hypergraph Tur\'an problems of higher uniformity.

\begin{thm}\label{nonprincip}
    For every integer $k\ge 4$ there exists a pair $P_k,P_k'$ of posets such that
    \begin{enumerate}
        \item 
    $e(P_k)=e^*(P_k)=e(P'_k)=e^*(P_k')=k$, 
    \item 
    $\pi^-(P_k),\pi^{*,-}(P_k),\pi^-(P'_k),\pi^{*,-}(P'_k)>k$,
    \item 
    and $\La(n,\{P_k,P'_k\})=\sum_{i=0}^{k-1}\binom{n}{\lfloor \frac{n-k}{2}\rfloor+i}$ and thus $\pi(\{P_k,P'_k\})=\pi^*(\{P_k,P'_k\})=k$.
    \end{enumerate}
\end{thm}

The construction showing (2) is an alteration of that of Tompkins \cite{T}. The upper bound proof of (3) uses a result of Tompkins and Wang \cite{TW} and Martin, Methuku, Uzzell, and Walker \cite{MMUW} and the chain-partition method \cite{GLL,GL}.

\section{Proof}
We start with some definitions. The diamond poset $\Diamond$ has four elements $a,b,c,d$ with $a<b<d,a<c<d$ being all its cover relations. Let $\vee$ be the poset on three elements $a,b,c$ with $a<b,c$ being all its strict relations, and $\wedge$ be its dual, i.e. the poset on $a,b,c$ with $a>b,c$.

For two posets $(P,\leqslant_P)$ and $(Q,\leqslant_Q)$, let their \textit{ordinal sum} $P + Q$ be the poset on the disjoint union of $P$ and $Q$ with $p\leqslant_{P+ Q} q$ for all $p\in P,q\in Q$ and the relations within $P$ and $Q$ inherited from $\leqslant_P$ and $\leqslant_Q$, respectively.

Let $Y_k$ be the poset $C_{k-1}+\vee$ and let $Y'_k$ be the poset $\wedge+C_{k-1}$, where $C_k$ is the chain on $k$ elements.

Finally, for $\F\subseteq 2^{X}$ with $|X|=n$ let $\lambda_n(\F)=\sum_{F\in \F}\binom{n}{|F|}^{-1}$ be the \textit{Lubell mass} of $\F$. Note that $\lambda_n(\F)=\mathbb{E}(|\C\cap \F|)$, where $\C$ is a random maximal chain chosen uniformly from $\bC_X$ the set of all maximal chains in $2^{X}$. Also, $\lambda_n(\F)\le t$ implies $|\F|\le t\binom{n}{\lfloor \frac{n}{2}\rfloor}$. 

\begin{proof}[Proof of Theorem \ref{nonprincip}]
    For $k\ge 4$, let us define $P_k=C_{k-2}+\Diamond=C_1+Y_{k-2}+C_1$ and $P'_k=\Diamond+C_{k-2}=C_1+Y'_{k-2}+C_1$.

    First we claim $e(P_k)=e^*(P_k)=e(P'_k)=e^*(P'_k)=k$. As $P_k,P'_k$ contain chains of length $k+1$, we have $e^*(P_k)\ge e(P_k)\ge k$ and $e^*(P'_k)\ge e(P'_k)\ge k$. On the other hand, $k+1$ consecutive layers contain a saturated chain $F_0\subsetneq F_1\subsetneq\dots \subsetneq F_k$ of length $k+1$ and together with the set $G_{k-1}\neq F_{k-1}$, $F_{k-2}\subsetneq G_{k-1}\subsetneq F_k$ they form a strong copy of $P_k$ and together with the set $G_1\neq F_1$, $F_0\subsetneq G_1\subsetneq F_2$ they form a strong copy of $P'_k$. This shows $e^*(P_k),e^*(P'_k)\le k$.

    Tompkins introduced \cite{T} a $\Diamond$-free family $\F\subseteq \binom{[n]}{\lfloor \frac{n}{2}\rfloor-1}\cup \binom{[n]}{\lfloor \frac{n}{2}\rfloor}\cup \binom{[n]}{\lfloor \frac{n}{2}\rfloor+1}$ of size $(2+\frac{29}{135}\Delta_4+o(1))\binom{n}{\lfloor \frac{n}{2}\rfloor}$, where $\Delta_q=\prod_{i=1}^\infty(1-q^{-i})$. Let $\F_k=\F\cup_{i=-k+1}^{-2} \binom{[n]}{\lfloor\frac{n}{2}\rfloor+i}$. Clearly, $|\F_k|=(k+\frac{29}{135}\Delta_4+o(1))\binom{n}{\lfloor \frac{n}{2}\rfloor}$. We claim that $\F_k$ is weak $P_k$-free. Indeed, as $\F_k$ lives on $k+1$ levels, each maximal chain in a copy of $P_k$ in $\F_k$ must contain one set from each level involved. Thus the $\Diamond$ part of $P_k=C_{k-2}+\Diamond$ must lie in $\F$, but that is $\Diamond$-free. Similarly, $\bar{\F_k}=\{[n]\setminus F:F\in\F_k\}$ is weak $P'_k$-free. This proves $\pi^-(P_k),\pi^{*,-}(P_k),\pi^-(P'_k),\pi^{*,-}(P'_k)>k$.

    Finally, let $\G\subseteq 2^{[n]}$ be a strong $\{P_k,P'_k\}$-free family. Let us use the min-max partition $$\bC_n=\bC_-\cup \bigcup_{G\in \G}\bC_G\cup\bigcup_{G_1\subsetneq G_2, ~G_1,G_2\in\G}\bC_{G_1,G_2},$$ where $\bC_-$ is the set of maximal chains $\C$ with $\C\cap \G=\emptyset$, $\bC_G$ is the set of maximal chains with $\C\cap \G=\{G\}$, and $\bC_{G_1,G_2}$ is the set of maximal chains with $G_1$ being the smallest, and $G_2$ being the largest set in $\C\cap \G$. As $P_k=C_1+Y_{k-2}+C_1,P'_k=C_1+Y'_{k-2}+C_1$ and $\G$ is strong $\{P_k,P'_k\}$-free, for any $G_1\subsetneq G_2, ~G_1,G_2\in \G$, the family $\G^0_{G_1,G_2}=\{G\setminus G_1:G\in \G, ~G_1\subsetneq G\subsetneq G_2\}$ is strong $\{Y_{k-2},Y'_{k-2}\}$-free. Theorem 4 of the papers \cite{MMUW,TW} states that the Lubell mass of a strong $\{Y_\ell,Y'_\ell\}$-free family $\HH \subseteq 2^X$ with $\emptyset,X\notin \HH$ is at most $\ell$ provided $|X|\ge \ell+1$ and $\ell\ge 2$. If $|G_2|-|G_1|\ge k$, then applying this result to $\G^0_{G_1,G_2}$ with $\ell=k-2$, we obtain $\lambda_{|G_2|-|G_1|}(\G^0_{G_1,G_2})\le k-2$ and therefore $\lambda_{|G_2|-|G_1|}(\G_{G_1,G_2})\le k$ where $\G_{G_1,G_2}=\{G\setminus G_1:G\in \G, ~G_1\subseteq G\subseteq G_2\}$. If $|G_2|-|G_1|\le k-1$, then $\lambda_{|G_2|-|G_1|}(\G_{G_1,G_2})\le k$ is trivially true. 
    
    Conditional on belonging to a nonempty class
$\bC_{G_1,G_2}$, the segment of $\C$ between $G_1$ and $G_2$
is uniformly distributed among the maximal chains of that interval:
the restrictions below $G_1$ and above $G_2$ do not depend on
the internal segment. As a consequence, $\lambda_n(\G)=\mathbb{E}(|\C\cap \G|)$ is a weighted average of the $\lambda_{|G_2|-|G_1|}(\G_{G_1,G_2})$s and 0 (for $\bC_-$) and 1 (for $\bC_G$s). Therefore, $\lambda_n(\G)\le k$ and so (see \cite[Lemma 3.2]{GLL}) $|\G|\le \sum_{i=0}^{k-1}\binom{n}{\lfloor \frac{n-k}{2}\rfloor+i}$. Conversely, the union of the $k$ middle levels contains no chain
on $k+1$ elements, so it is weak $\{P_k,P'_k\}$-free and has size
$\sum_{i=0}^{k-1}\binom{n}{\lfloor \frac{n-k}{2}\rfloor+i}$. This finishes the proof of (3).
\end{proof}

\end{document}